\documentclass[12pt]{amsart}

\makeatletter
\@namedef{subjclassname@2020}{\textup{2020} Mathematics Subject Classification}
\makeatother

\newcommand{\pii}{\pi i}

\newcommand{\dsum}{\di\sum}

\DeclareFontFamily{U}{mathx}{\hyphenchar\font45}
\DeclareFontShape{U}{mathx}{m}{n}{
      <5> <6> <7> <8> <9> <10>
      <10.95> <12> <14.4> <17.28> <20.74> <24.88>
      mathx10
      }{}
\DeclareSymbolFont{mathx}{U}{mathx}{m}{n}
\DeclareFontSubstitution{U}{mathx}{m}{n}
\DeclareMathAccent{\widecheck}{0}{mathx}{"71}

\newcommand{\inrr}{\ensuremath{\in\rr}}

\renewcommand{\kill}[1]{}
\newcommand{\dummy}[1]{\mbox{}}

\makeatletter
\newcommand{\xequal}[2][]{\ext@arrow 0055{\equalfill@}{#1}{#2}}
\def\equalfill@{\arrowfill@\Relbar\Relbar\Relbar}
\makeatother

\renewcommand{\k}{\ensuremath{\ol{\mathrm{P}}}}

\newcommand{\h}{\hline}

\renewcommand{\k}[1]{\ensuremath{\left({#1}\right)}}

\newcommand{\re}{\te{Re}\,}

\newcommand{\bca}{\begin{cases}}
\newcommand{\eca}{\end{cases}}

\newcommand{\mug}{\ensuremath{\infty}}

\newcommand{\sinx}{\ensuremath{\sin x}}

\newcommand{\ff}[2]{\ensuremath{\di\fr{#1}{#2}}}

\renewcommand{\ss}[3]{\ensuremath{\di\int_{#1}^{#2}{#3}\,dx}}

\newcommand{\bpic}{\begin{picture}}\newcommand{\epic}{\end{picture}}

\newcommand{\beda}{\begin{edaenumerate}}
\newcommand{\eeda}{\end{edaenumerate}}

\newcommand{\cd}{\cdots}

\newcommand{\asx}{\ensuremath{\sin^{-1}x}}

\newcommand{\acx}{\ensuremath{\cos^{-1}x}}

\newcommand{\q}{\quad}

\newcommand{\bq}{\begin{quote}}\newcommand{\eq}{\end{quote}}

\newcommand{\rt}{\sqrt}
\newcommand{\be}{\begin{enumerate}}\newcommand{\ee}{\end{enumerate}}
\newcommand{\bce}{\begin{center}}\newcommand{\ece}{\end{center}}
\newcommand{\bde}{\begin{description}}\newcommand{\ede}{\end{description}}
\newcommand{\bri}{\begin{flushright}}\newcommand{\eri}{\end{flushright}}
\newcommand{\bb}{\begin{block}}\newcommand{\eb}{\end{block}}
\newcommand{\bt}{\begin{thm}}\newcommand{\et}{\end{thm}}
\newcommand{\bpf}{\begin{proof}}\newcommand{\epf}{\end{proof}}
\newcommand{\bex}{\begin{ex}}\newcommand{\eex}{\end{ex}}
\newcommand{\bexr}{\begin{exr}}\newcommand{\eexr}{\end{exr}}
\newcommand{\bft}{\begin{fact}}\newcommand{\eft}{\end{fact}}
\newcommand{\brk}{\begin{rmk}}\newcommand{\erk}{\end{rmk}}
\newcommand{\ba}{\begin{align*}}\newcommand{\ea}{\end{align*}}
\newcommand{\bexe}{\begin{exe}}\newcommand{\eexe}{\end{exe}}

\newcommand{\bit}{\begin{itemize}}\newcommand{\eit}{\end{itemize}}

\newcommand{\bcm}{}

\newcommand{\fr}{\frac}
\newcommand{\cc}{\ensuremath{\mathbf{C}}}

\newcommand{\rr}{\ensuremath{\mathbf{R}}}

\newcommand{\bd}{\begin{defn}}\newcommand{\ed}{\end{defn}}
\newcommand{\bp}{\begin{prop}}\newcommand{\ep}{\end{prop}}
\newcommand{\p}{\ensuremath{\pi}}
\newcommand{\eh}{\emph}

\newcommand{\te}{\text}

\newcommand{\di}{\displaystyle}

\newcommand{\f}{\frac}

\newcommand{\z}{\ensuremath{\bm{z}}}

\usepackage[dvipdfmx]{graphicx}
\usepackage[dvipsnames]{xcolor}
\usepackage{float,afterpage}
\usepackage{asymptote,layout}
\usepackage{wrapfig,epic}

\usepackage{geometry}
\usepackage{exscale,latexsym,bm}
\usepackage{amssymb,enumerate,amsmath,amsthm,amsfonts}
\usepackage{verbatim,fancybox,wasysym,fancyhdr,type1cm}
\usepackage[frame,all,poly,curve,knot,arrow]{xy}
\usepackage{colortbl}
\usepackage{boxedminipage}
\usepackage{multirow}

\graphicspath{{./figs/}}
\theoremstyle{definition}
\newtheorem{thm}{Theorem}[section]

\newtheorem*{mth}{Main Theorem}

\newtheorem{prop}[thm]{Proposition}\newtheorem{cor}[thm]{Corollary}

\newtheorem{exr}[thm]{Exercise}

\newtheorem{ex}[thm]{Example}

\newtheorem{defn}[thm]{Definition}\newtheorem{rmk}[thm]{Remark}
\newtheorem{fact}[thm]{Fact}
\newtheorem{block}[thm]{}
\newtheorem*{exe}{Exercise}

\renewcommand{\z}{\zeta}

\begin{document}
\renewcommand{\h}{\hline}
\renewcommand{\arraystretch}{1.5}
\newcommand{\fib}{\text{fib\,}}
\newcommand{\fdp}{\text{FDP\,}}

\title{Integral representations for $\zeta(3)$ with the inverse sine function}
\author{Masato Kobayashi}
\date{\today}                                       
\address{Masato Kobayashi\\
Department of Engineering\\
Kanagawa University, 3-27-1 Rokkaku-bashi, Yokohama 221-8686, Japan.}
\email{masato210@gmail.com}


\maketitle
\begin{abstract}
We show four new integral representations for $\zeta(3)$ as a reformulation of Ewell (1990) and Yue-Williams (1993) with the inverse sine function and Wallis integral.
As a consequence, we also show a local integral representation for the trilogarithm function. 
\end{abstract}
\tableofcontents

\renewcommand{\asx}{\sin^{-1}x}
\renewcommand{\acx}{\cos^{-1}x}

\section{Introduction}

\subsection{Ap\'{e}ry number}
The Riemann zeta function is one of the important topics in number theory. 
For a complex number $s$ such that 
$\re{(s)}>1$, we usually define 
\[
\z(s)=
\ff{1}{1^{s}}+\ff{1}{2^{s}}+\ff{1}{3^{s}}+\ff{1}{4^{s}}+\cd.
\]
Leonhard Euler proved 
\[
\z(2n)=-\ff{1}{2}\ff{(2\pii)^{2n}}{(2n)!}B_{2n} \q (n\ge1)
\]
($\{B_{2n}\}$ are Bernoulli numbers) while 
$\{\z(2n+1)\}_{n\ge 1}$ remain to be unknown. 
In this article, we particularly study 
\[
\zeta (3)=1.2020569\cdots,
\]
the \eh{Apery number}, as Ap\'{e}ry proved its irrationality \cite{apery} in 1979. 
Although its exact value is still unclear, 
there exist many representations for $\z(3)$ such as 
\begin{align*}
	\z(3)&=\ff{2}{7}\p^{2}\log2+\ff{16}{7}
	\di\int_{0}^{\p/{2}}{x\log (\sinx)}\,dx,\\
	\z(3)&=\ff{7}{180}\p^{3}-
2\dsum_{n=1}^{\mug} \ff{1}{n^{3}(e^{2\pi n}-1)},
	\\
	\z(3)&=\ff{5}{2}
	\dsum_{n=1}^{\mug}\ff{(-1)^{n-1}}{n^{3}\binom{2n}{n}}
\end{align*}
due to Euler, Plouffe \cite{plouffe}, Ramanujan, Ap\'{e}ry and doubtless others; see also Chen-Srivastava \cite{chen} and  Nash-O'Connor \cite{nash} for other representations.

\subsection{Main result}
The main result of this article is to show 
the following integral representations for $\z(3)$; 
all of these are new. 
\begin{mth}\label{mth}
\begin{align*}
	\z(3)&=\ff{8}{7}\di\int_{0}^{1}\ff{\sin^{-1}x\cos^{-1}x}{x}dx\\
	&=\ff{8}{\p}\di\int_{0}^{1}\ff{(\sin^{-1}x)^{2}\cos^{-1}x}{x}dx\\
	&=\ff{16}{5\p}\di\int_{0}^{1}\ff{\sin^{-1}x(\cos^{-1}x)^{2}}{x}dx\\
	&=\ff{32}{3\p}\di\int_{0}^{1}\ff{(\sinh^{-1}x)^{2}\cos^{-1}x}{x}dx.\\
\end{align*}
\end{mth}
Before going into the proof of this theorem in the next section, we wish to mention work of Ewell (1990) and Yue-Williams (1993) and set up notation.

\subsection{Results of Ewell and Yue-Williams}

\begin{fact}
[Ewell \cite{ewell}, Yue-Williams \te{\cite[p.1582]{williams}}]
\[
\z(3)=\ff{\p^{2}}{7}
\k{1-4
\dsum_{n=1}^{\mug}\ff{\z(2n)}{(2n+1)(2n+2)2^{2n}}
},
\]
\[
\z(3)=-2\p^{2}
\k{\dsum_{n=0}^{\mug}\ff{\z(2n)}{(2n+2)(2n+3)2^{2n}}
}.
\]
\end{fact}
These series are quite similar because they both derived some infinite sums related to $\z(3)$ 
from the Maclaurin series involving $\asx$ using 
Wallis integral; this method is modification of Boo Rim Choe's elementary proof for $
\sum_{n=1}^{\mug}\f{1}{n^{2}}
=\f{\p^{2}}{6}$ \cite{boo}. 


\subsection{Notation}

Throughout $n$ denotes a nonnegative integer.
Let 
\begin{align*}
	(2n)!!&=2n(2n-2)\cd 4\cdot 2,
	\\(2n-1)!!&=(2n-1)(2n-3)\cd 3\cdot 1.
\end{align*}
In particular, we understand that $(-1)!!=0!!=1$.
Moreover, let \[
c_{n}=\ff{(n-1)!!}{n!!}.
\]
This number appears in the following integral:
\begin{fact}[Wallis integral]
\[
\ss{0}{\pi/2}{\sin^n x}
=
\begin{cases}
	\f{\p}{2} c_{n}&	\text{$n$ even,}\\
	c_{n}&	\text{$n$ odd.}\\
\end{cases}
\]
\end{fact}
By $\sin^{-1} x$ and $\cos^{-1} x$, we mean the real inverse sine and cosine functions ($\arcsin x, \arccos x$), that is, 
\[
\begin{array}{ccl}
	y=\asx&\iff   &   x=\sin y, \q-\f{\,\pi\,}{2}\le y\le \f{\,\pi\,}{2},\\
	y=\acx&\iff   &    x=\cos y, \q0\le y\le \p.
\end{array}
\]
\begin{fact}
\[
\asx=
\dsum_{n=0}^{\mug}c_{2n}\ff{x^{2n+1}}{2n+1}, \q |x|\le 1.
\]\end{fact}
\begin{rmk}
In the literature, some researchers exclude $|x|=\pm1$. 
However, even for $|x|=\pm1$, this equality indeed holds as Boo Rim Choe pointed out \cite{boo}.
\end{rmk}
Further, 
$\sinh^{-1}x=\log(x+\rt{x^{2}+1})$ $(x\inrr)$ 
denotes the inverse hyperbolic sine function 
(some authors write $\te{arsinh\,} x,\te{arcsinh\,} x$ or $\te{argsinh\,} x$ for this function).

\section{Proof of Main theorem}
In this section, we give a proof of 
equalities in Main Theorem one by one.

\begin{thm}\label{th1}
\[
\z(3)=
\ff{8}{7}\di\int_{0}^{1}\ff{\sin^{-1}x\cos^{-1}x}{x}dx.
\]
\end{thm}

\begin{proof}
Let us start with the Maclaurin series
\[
\ff{\sin^{-1}y}{y}=
\dsum_{n=0}^{\mug}c_{2n}\ff{y^{2n}}{2n+1}, \q |y|\le 1.
\]
Let $t, u$ be real variables such that 
$0\le t, u\le 1$. 
Recall from calculus that term-wise integration is possible for convergent power series. By integrating the series above 
from 0 to $tu$, we have 
\begin{align*}
	I(t, u):&=\int_{0}^{tu}\ff{\sin^{-1}y}{y}dy
=
\int_{0}^{tu}
\k{
\dsum_{n=0}^{\mug}c_{2n}\ff{y^{2n}}{2n+1}}
dy	
	\\&=\dsum_{n=0}^{\mug}c_{2n}\ff{1}{2n+1}
\int_{0}^{tu}y^{2n}dy
=
\dsum_{n=0}^{\mug}c_{2n}\ff{1}{(2n+1)^{2}}{t^{2n+1}u^{2n+1}}{}.
\end{align*}
Now, let $v=\sin^{-1}u$ $(0\le v\le \f{\,\pi\,}{2})$.
Then $u=\sin v$.
Integrate $I(t, \sin v)$ from 0 to $\f{\,\pi\,}{2}$ 
in $v$. On one hand, 
\begin{align*}
	J(t):&=\int_{0}^{\p/2}I(t, \sin v)dv
	=
\dsum_{n=0}^{\mug}c_{2n}\ff{t^{2n+1}}{(2n+1)^{2}}
\underbrace{\int_{0}^{\pi/2}{\sin^{2n+1} v}dv}_{
c_{2n+1}=\f{1}{c_{2n}(2n+1)}
}
	\\&=
\dsum_{n=0}^{\mug}\ff{t^{2n+1}}{(2n+1)^{3}}.
\end{align*}
On the other hand, exchanging order of the integrals yields
\begin{align*}
	J(t)&=\int_{0}^{\p/2}
\int_{0}^{t\sin v}
\ff{\sin^{-1}y}{y}dydv
=
\int_{0}^{t}
\int_{\sin^{-1}\f{y}{t}}^{\p/2}
\ff{\sin^{-1}y}{y}dvdy
	\\&=\int_{0}^{t}
\ff{\sin^{-1}y}{y}
\k{ \ff{\,\pi\,}{2}-\sin^{-1}\f{y}{t}}
dy
=
\int_{0}^{1}
\ff{\sin^{-1}(tx)\cos^{-1}x}{x}dx.
\end{align*}
Thus, 
\[
(1-2^{-3})\z(3)
=\dsum_{n=0}^{\mug}\ff{1}{(2n+1)^{3}}
=J(1)=
\int_{0}^{1}
\ff{\sin^{-1}x\cos^{-1}x}{x}dx.
\]
Conclude that 
\[
\z(3)=
\ff{8}{7}\di\int_{0}^{1}\ff{\sin^{-1}x\cos^{-1}x}{x}dx.
\]
\end{proof}

\begin{cor}
\[
\di\int_{0}^{1}\ff{\sinh^{-1}x\acx}{x}dx=
\ff{\p^{3}}{32}.
\]
\end{cor}
\begin{proof}
Start with
\[
	\sinh^{-1}y=
	\dsum_{n=0}^{\mug}(-1)^{n}
	c_{2n}\ff{y^{2n+1}}{2n+1}, \q |y|\le 1
\]
which easily follows from 
\[
\sin^{-1}y=
\dsum_{n=0}^{\mug}
c_{2n}\ff{y^{2n+1}}{2n+1}, \q |y|\le 1
\]
and $\sinh^{-1}z=-i\sin^{-1}(iz)$ (for all $z\in\cc$) 
\cite[p.87]{ab}.
Now consider the argument with replacing $\sin^{-1}y$ by $\sinh^{-1}y$ in Proof of Theorem \ref{th1} 
throughout. 
Then, the proof goes without any substantial changes.
 It leads us to 
\[
\dsum_{n=0}^{\mug}\ff{(-1)^{n}}{(2n+1)^{3}}=
\di\int_{0}^{1}\ff{\sinh^{-1}x\acx}{x}dx.
\]
The left hand side is $\f{\p^{3}}{32}$ 
\cite[p.808]{ab}. 
\end{proof}

\begin{thm}\label{th2}
\[
\z(3)=
\ff{8}{\p}
\int_{0}^{1}
\ff{(\sin^{-1}x)^{2}\cos^{-1}x}{x}dx.
\]
\end{thm}

\begin{proof}
Start with 
\[
\ff{(\sin^{-1}y)^{2}}{y}=
\ff{\,1\,}{2} \dsum_{n=1}^{\mug}\ff{1}{c_{2n}}\ff{y^{2n-1}}{n^{2}}, \q |y|\le 1.
\]
Again, let $t, u$ be real variables such that $0\le t, u\le 1$. Set 
\[
K(t, u):=\di\int_{0}^{tu}\ff{(\sin^{-1}y)^{2}}{y}dy.
\]
Then 
\begin{align*}
	K(t, u)&=\di\int_{0}^{tu}
\k{
\ff{\,1\,}{2} \dsum_{n=1}^{\mug}\ff{1}{c_{2n}}\ff{y^{2n-1}}{n^{2}}}
dy
=
\ff{\,1\,}{2} \dsum_{n=1}^{\mug}\ff{1}{c_{2n}}\ff{1}{n^{2}}
\di\int_{0}^{tu}{y^{2n-1}}\,dy
	\\&=\ff{\,1\,}{4}
\dsum_{n=1}^{\mug}\ff{1}{c_{2n}}\ff{t^{2n}u^{2n}}{n^{3}}.
\end{align*}
Let $v=\sin^{-1}u$ $(0\le v\le \f{\,\pi\,}{2})$
. Then $u=\sin v$.
Integrate $K(t, \sin v)$ from 0 to $\f{\,\pi\,}{2}$ in $v$:
\begin{align*}
	L(t):&=\int_{0}^{\p/2}K(t, \sin v)dv=
\ff{\,1\,}{4} \dsum_{n=1}^{\mug}\ff{1}{c_{2n}}\ff{t^{2n}}{n^{3}}
\underbrace{\int_{0}^{\pi/2}{\sin^{2n} v}dv}_{c_{2n}\f{\p}{2}}
	\\&=\ff{\p}{8}
\dsum_{n=1}^{\mug}\ff{t^{2n}}{n^{3}}.
\end{align*}
On the other hand, 
exchanging order of the integrals yields
\begin{align*}
	L(t)&=\int_{0}^{\p/2}
\int_{0}^{t\sin v}
\ff{(\sin^{-1}y)^{2}}{y}dydv
=
\int_{0}^{t}
\int_{\sin^{-1}\f{y}{t}}^{\p/2}
\ff{(\sin^{-1}y)^{2}}{y}dvdy
	\\&=\int_{0}^{t}
\ff{(\sin^{-1}y)^{2}}{y}
\k{
\ff{\,\pi\,}{2}-
\sin^{-1}\ff{y}{t}
}
dy
	\\&=\int_{0}^{1}
\ff{(\sin^{-1}(tx))^{2}\cos^{-1}x}{x}
dx.
\end{align*}
Hence 
\[
\ff{\p}{8}\z(3)=
\ff{\p}{8}
\dsum_{n=1}^{\mug}\ff{1}{n^{3}}=
L(1)=
\int_{0}^{1}
\ff{(\sin^{-1}x)^{2}\cos^{-1}x}{x}dx
\]
and conclude that 
\[
\z(3)=
\ff{8}{\p}
\int_{0}^{1}
\ff{(\sin^{-1}x)^{2}\cos^{-1}x}{x}dx.
\]
\end{proof}


%

\begin{thm}\label{th4}
\[
\z(3)=\ff{32}{3\p}
\di\int_{0}^{1}\ff{(\sinh^{-1}x)^{2}\acx}{x}dx.
\]
\end{thm}

\begin{proof}
With $\sinh^{-1}z=-i\sin^{-1}(iz)$ (for $z\in\cc$), we see that 
\[
\ff{(\sinh^{-1}y)^{2}}{y}=
\ff{1}{2}
\dsum_{n=1}^{\mug}\ff{(-1)^{n-1}}{c_{2n}}\ff{y^{2n-1}}{n^{2}}, \q |y|\le 1.
\]
Again, consider the argument with replacing $\sin^{-1}y$ by $\sinh^{-1}y$ in Proof of Theorem \ref{th2}  throughout. 
Then, the proof goes without any substantial changes.
It leads us to 
\[
\ff{\p}{8}
\dsum_{n=1}^{\mug}\ff{(-1)^{n-1}}{n^{3}}=
\di\int_{0}^{1}\ff{(\sinh^{-1}x)^{2}\acx}{x}dx,
\]
\[
\ff{\p}{8}
(1-2^{1-3})\z(3)=
\di\int_{0}^{1}\ff{(\sinh^{-1}x)^{2}\acx}{x}dx,
\]
and hence 
\[
\z(3)=\ff{32}{3\p}
\di\int_{0}^{1}\ff{(\sinh^{-1}x)^{2}\acx}{x}dx.
\]
\end{proof}


\begin{cor}\label{th5}
\[
\z(3)=
\ff{16}{5\p}
\di\int_{0}^{1}\ff{\asx(\acx)^{2}}{x}dx.
\]
\end{cor}

\begin{proof}
This is an easy consequence of Theorems \ref{th1} and \ref{th2} as follows. 
Recall that $\asx+\acx=\f{\p}{2}$ whenever $0\le x\le 1$. Thus, 
\begin{align*}
	\di\int_{0}^{1}\ff{\asx(\acx)^{2}}{x}dx&=
	\di\int_{0}^{1}\ff{\asx\acx}{x}
\k{ \ff{\,\pi\,}{2}-\asx}
dx
	\\&=\ff{\,\pi\,}{2}
\di\int_{0}^{1}\ff{\asx\acx}{x}dx
-
\di\int_{0}^{1}\ff{(\asx)^{2}\acx}{x}dx
	\\&=\ff{\,\pi\,}{2}\k{\ff{7}{8}\z(3)}-\ff{\p}{8}\z(3)
=
\ff{5\p}{16}\z(3).
\end{align*}
Hence 
\[
\z(3)=
\ff{16}{5\p}
\di\int_{0}^{1}\ff{\asx(\acx)^{2}}{x}dx.
\]
\end{proof}
We completed the proof of Main Theorem.

\section{Further remarks}

We end with some remarks for our future research.
\begin{enumerate}
	\item Main Theorem suggests us to think of the following family of integrals:
\[
I(m, n)=
\di\int_{0}^{1}{\ff{(\asx)^{m}(\acx)^{n}}{x}}\,dx, \q 
m, n\ge 0, m+n\ge 1.
\]
Here, let us observe several examples for some interest. 
As is well-known, 
\[
I(1, 0)=
\di\int_{0}^{1}{\ff{\asx}{x}}\,dx=\ff{1}{2}\pi\log2
\]
while 
\[
I(0, 1)=\di\int_{0}^{1}{\ff{\acx}{x}}\,dx
=\di\int_{0}^{1}\k{
\ff{\pi}{2}\ff{1}{x}-1-\ff{x^{2}}{6}-\cd
}\,dx
\]
is divergent. 
For $m, n\ge 1$, the first three integrals are 
\[
I(1, 1)=\ff{7}{8}\z(3), \q
I(2, 1)=\ff{\p}{8}\z(3), \q
I(1, 2)=\ff{5\p}{16}\z(3)
\]
as we have shown before. 
Now such integrals satisfy the relations
\[
I(m, n)=\ff{\p}{2}I(m-1, n)-I(m-1, n+1), \q m\ge 1
\]
and similarly 
\[
I(m, n)=\ff{\p}{2}I(m, n-1)-I(m+1, n-1), \q n\ge 1
\]
since $\asx+\acx= \f{\,\pi\,}{2}$ for $0\le x\le 1$. 
We then have 
\[
I(2, 0)=\ff{\p}{2}I(1, 0)-I(1, 1)=
\ff{\p^{2}}{4}\log2-\ff{7}{8}\z(3),
\]
\[
I(3, 0)=\ff{\p}{2}I(2, 0)-I(2, 1)=
\ff{\p^{3}}{8}\log2-\ff{9}{16}\z(3)
\]
and so on.
Wolfram alpha \cite{wolfram} says that 
$I(4, 0)$ ``is close to"
\[
\ff{1}{32}\k{-18\p^{2}\z(3)+93\z(5)+2\p^{4}\log2}
\]
which involves $\z(3)$ and $\z(5)$; 
actually, we computed this as 
\[
I(4, 0)=
\di\int_{0}^{1}{\ff{(\asx)^{4}}{x}}\,dx
=
\di\int_{0}^{\p/2}{u^{4}\cot u}\,du.
\]
Once we can figure out all the coefficients of $\f{(\asx)^{4}}{x}$, then it may be possible to 
derive some series involving $\z(3)=\sum_{n=1}^{\mug}\f{1}{n^{3}}$ and 
$\z(5)=\textstyle\sum_{n=1}^{\mug}\f{1}{n^{5}}$ in a  similar method.
	\item The \emph{polylogarithm function} $\te{Li}_{s}(t)=
\sum_{n=1}^{\mug}\f{t^{n}}{n^{s}}$ plays a significant role in many areas of number theory.
	As a byproduct of Proof of Theorem \ref{th2}, we obtained a local integral representation of the trilogarithm function; 
that is, for $0\le t\le 1$, we have 
\[
\te{Li}_{3}(t)
=
\ff{8}{\p}
\di\int_{0}^{1}\ff{(\sin^{-1}(\rt{t}x))^{2}\acx}{x}dx.
\]
This seems to be also new.
\end{enumerate}

\begin{center}
\textbf{Acknowledgment}
\end{center}

\begin{center}
This research arose from Iitaka online seminar in 2020-2021. 
The author would like to thank the organizer 
Shigeru Iitaka, Yuji Yamaga and Kouichi Nakagawa for fruitful discussions. He also thanks Satomi Abe as well as 
Michihito Tobe for supporting his writing the manuscript.
\end{center}


\begin{thebibliography}{99}
\bibitem{ab} M. Abramowitz, A. Stegun eds. (1972), Handbook of Mathematical Functions with Formulas, Graphs, and Mathematical Tables, New York: Dover Publications.
\bibitem{apery} R. Ap\'{e}ry, ``Irrationalit\'{e} de $\zeta(2)$ et $\zeta(3)$", Ast\'{e}risque {\textbf {61}} (1979), 11-13.
\bibitem{boo} Boo Rim Choe, An elementary proof of 
$\sum_{n=1}^{\mug} \f{1}{n^{2}}=\f{\p^{2}}{6}$, Amer. Math. Monthly \textbf{94} (1987), 662-663.
\bibitem{chen} M.-P. Chen and H. M. Srivastava, Some families of series representations for the Riemann $\z(3)$, Results Math. \textbf{33} (1998), 179-197.
\bibitem{ewell} J. A. Ewell, A new series representation for $\z(3)$, Amer. Math. Monthly {\textbf {97}} (1990), 219-220. 
\bibitem{nash} C. Nash and D. O'Connor, Ray-Singer torsion, topological field theories and the Riemann Zeta function at $s=3$, in Low-Dimensional Topology and Quantum Field Theory, Plenum Press, New York and London, 1993, 279-288. 
\bibitem{plouffe} S. Plouffe, Identities inspired from Ramanujan Notebooks II, 1998.
\bibitem{williams} K. S. Williams and Z.N. Yue, Some series representations of $\z(2n+1)$, Rocky Mt. J. Math.  
\textbf{23} (1993), 1581-1592.
\bibitem{wolfram} Wolfram alpha, 
\texttt{https://www.wolframalpha.com}.
\end{thebibliography}
\end{document}